\documentclass[11pt]{amsart}

\usepackage{times}
\usepackage{amsmath,  amssymb,slashed,url,bm,upgreek,amsthm,esint}
\usepackage{graphicx,enumerate}
\usepackage{hyperref}
\newcounter{intro}

\newtheorem{theo}[intro]{Theorem}
\newtheorem{coro}[intro]{Corollary}
\newtheorem{propo}[intro]{Proposition}

\newtheorem{thm}{Theorem}[section]
\newtheorem{lem}[thm]{Lemma}
\newtheorem{prop}[thm]{Proposition}
\newtheorem{cor}[thm]{Corollary}

\newtheorem{rems}[thm]{Remarks}

\newcommand{\cref}[1]{Corollary~\ref{#1}}

\newcommand{\pref}[1]{Proposition~\ref{#1}}

\newcommand{\tref}[1]{Theorem~\ref{#1}}

\DeclareMathOperator{\Id}{Id}

\DeclareMathOperator{\trace}{Trace}
\DeclareMathOperator{\spec}{spec}

\DeclareMathOperator{\ricci}{Ricci}
\DeclareMathOperator{\ricm}{\rho_-}

\DeclareMathOperator{\diam}{diam}

\DeclareMathOperator{\dv}{dv}

\DeclareMathOperator{\vol}{vol}\DeclareMathOperator{\dvol}{dvol}
\def\R{\mathbb R}\def\N{\mathbb N}\def\bS{\mathbb S}\def\bB{\mathbb B}

\def\cC{\mathcal C}

\def\sch{Schr\"odinger }

\def\cA{\mathcal A}
\newcommand{\euler}{\mathrm{e}}
\newcommand{\drm}{\mathrm{d}}
\newcommand{\RR}{\mathbb{R}}
\newcommand{\NN}{\mathbb{N}}
\DeclareMathOperator{\Ric}{\mathop{Ric}}

\newcommand{\hmt}[1]{\leavevmode{\marginpar{\tiny%
 $ \hbox to 0mm{\hspace*{-0.5mm} $ \leftarrow $ \hss}%
 \vcenter{\vrule depth 0.1mm height 0.1mm width \the\marginparwidth}%
 \hbox to
 0mm{\hss $ \rightarrow $ \hspace*{-0.5mm}} $ \\\relax\raggedright #1}}}

\begin{document}
\title[Geometric and spectral estimates based on spectral {R}icci curvature assumptions]{Geometric and spectral estimates based on spectral {R}icci curvature assumptions}
\author{ Gilles Carron \and 
Christian Rose}
\email{Gilles.Carron@univ-nantes.fr}
\email{christian.rose@mathematik.tu-chemnitz.de}
\begin{abstract} We obtain a  Bonnet-Myers theorem under a spectral condition:   a closed Riemannian $(M^n,g)$ manifold for which the lowest eigenvalue of the Ricci tensor $\rho$ is such that the Schr\"odinger operator $(n-2)\Delta+\rho$ is positive has finite fundamental group. Further, as a continuation of our earlier results, we obtain isoperimetric inequalities from Kato-type conditions on the Ricci curvature. We also obtain the Kato condition for the Ricci curvature  under purely geometric assumptions.\end{abstract}
\date\today
\maketitle

\footnotetext{\emph{Mots cl\'es: } In\'egalit\'e de Sobolev, in\'egalit\'es isop\'erim\`etriques, }

%
%
\section{Introduction}
\subsection{} Several  classical results in differential geometry state geometric, spectral  or topological estimates for a closed Riemannian manifold  assuming pointwise curvature constraints. Prominent examples of such results depending on some lower bound on the Ricci curvature are the theorems of Bonnet-Myers, Bochner, estimates of isoperimetric constants, and eigenvalue estimates for the Laplace-Beltrami operator, as well as bounds on the Betti numbers. An intriguing question is whether the assumption of a sharp lower Ricci curvature bound can be relaxed to variable curvature bounds instead. Of course, one can ask why this should be interesting at all: the Ricci tensor as a continuous function on a compact set is bounded below anyway. The answer lies in the quantitative nature of our question: if we disturb the curvature assumptions of the classical theorems by a small portion of "bad" curvature, how do the classical estimates change? Within this context, the treatment of families of manifolds are of importance: One can ask for uniform properties of all of the members of the family. Indeed, one could make assumptions on the infimum of the Ricci curvatures of all the members and derive estimates, but they will be not stable under small changes and much weaker than perturbation statements. Results in this direction have been obtained by K-D. Elworthy and S. Rosenberg, \cite{ElworthyRosenberg-91}. From the pioneering work of S.~Gallot, there have been a lot of results based on integral pinching conditions for the Ricci curvature (\cite{Gallot-88,Gallot-88_2,Aubry-07,MR1487752,MR1487753,MR1684981,MR1714334,MR1742252, MR1709777, MR1775137, MR1934123, MR2191529, MR2505581}).\\
Recently, we have obtained heat kernel and eigenvalue estimates assuming only a Kato condition on the negative part of Ricci curvature \cite{Rose-17,Carron-16}. It has been shown that this type of condition is somewhat weaker that the integral pinching assumption made by S. Gallot. 
\subsection{}Let's us introduce several notations  and conventions: 
In all that follows $(M^n,g)$ is always a closed Riemannian manifold of dimension $n\in\NN$, i.e., compact and without boundary, and $\vol$, respectively $\cA$, denote the $n$-, respectively $(n-1)$-dimensional Hausdorff measure on $M$. Denote by $\Ric$ the Ricci tensor of $M$ considered pointwise as an symmetric endomorphism of $T^*M$. We let 
$$\rho\colon M\to\RR, x\mapsto \min\{\sigma(\Ric_x)\},$$
where $\sigma(A)$ denotes the spectrum of an operator $A$. We have $$\Ric_x\ge \rho(x) \Id_{T_xM}$$ or $\Ric_x\ge \rho(x) g_x$ in the sense of quadratic forms on the tangent space at $x\in M$.
If $x$ is a real number, its negative part will be denoted by $x_-=\max\{0, -x\}$, for instance:
 $$\Ric_x\ge -\rho_-(x) g_x\ \text{or}\ \Ric_x\ge \uplambda g-(\rho(x)-\uplambda)_- g_x.$$
As usual, we denote by $\Delta\geq 0$ the non-negative Laplace-Beltrami operator acting on functions defined on $M$. This operator has purely discrete spectrum 
$$\sigma(\Delta)=\{\lambda_i\mid i\in\NN_0\},$$
where $$0=\lambda_0<\lambda_1\leq\lambda_2\leq\ldots .$$
\subsection{}
Our first result is based on an a spectral assumption for a Schr\"odinger operator:
\begin{theo} Let $(M^n,g)$ be a closed Riemannian manifold of dimension $n\ge 3$.
\begin{enumerate}[i)]
\item For any $\upepsilon\in \left[0, \frac{3}{n+4}\right)$, there is an explicit constant $C(n,\upepsilon)>0$ satisfying $\lim_{\upepsilon\to 0+} C(n,,\upepsilon)=1$, such that if the Schr\"odinger operator
$$\upepsilon \Delta+\rho-(n-1)k^2$$ is non negative, then we have
$\diam (M,g)\le C(n,\upepsilon) \frac{\pi}{k}$.
\item If for some $\upepsilon\in \left[0,\frac{1}{n-2}\right]$, the bottom of the spectrum of Schr\"odinger operator
$$\upepsilon \Delta+\rho$$ is positive, then the fundamental group of $M$ is finite.
\end{enumerate}
\end{theo}

\begin{rems}
\begin{enumerate}
\item This result is a weak generalization of the Bonnet-Myers theorem that states that if a \textit{complete} Riemannian manifold $(M^n,g)$ satisfies
$$\Ric\ge (n-1)k^2,$$ then $M$ is closed, $\diam (M,g)\le \pi$ and the fundamental group of $M$ is finite.
\item With respect to integral bound on curvature, the sharp result has been obtained by E. Aubry (\cite{Aubry-07}):
Let $2p>n\geq 3$, there is an $\upepsilon(p,n)>0$ such that if a \textit{complete} Riemannian manifold $(M^n,g)$ satisfies
$$\frac 1{\vol(M)}\int_M(\rho-(n-1))_-^p \, \dvol\leq\upepsilon(p,n),$$
then $M$ is closed, there is an explicit estimate for the diameter of $(M,g)$ and the fundamental group of $M$ is finite.
\item The first result is obtained in two stages. One first proves a Sobolev inequality (using the idea of D. Bakry and I. Mondello (\cite{BakrySobolev, Mondello2017})), and then uses a result of D. Bakry and M. Ledoux (\cite{BakryLedouxDiam}).
\item For $n=3$, Bour and Carron (\cite{BourCarron2}) have shown that if the operator
$$2\Delta+\rho $$ is positive, then the universal cover of $M^3$ is $\bS^3.$
\item We do not know wether the threshold $1/(n-2)$ is sharp, but we believe it is not.
\item We do not know wether a complete Riemannian manifold $(M^n,g)$ such that for some positive $\upepsilon<1/(n-2)$, the Schr\"odinger operator
$$\upepsilon \Delta+\rho $$ is positive must be closed. However, the proof of the theorem would show that if a \textit{complete} Riemannian manifold $(M^n,g)$ of dimension $n\ge 3$ satisfies  a uniform lower bound on the Ricci curvature of the form $\Ric_x\ge -(n-1)k^2$, a non collapsing assumption $\inf_x \vol B(x,1)>0$, and the fact that for some $\upepsilon\in \left[0,\frac{1}{n-2}\right)$, the bottom of the spectrum of Schr\"odinger operator
$\upepsilon \Delta+\rho$ is positive, then $M$ is  closed.
\end{enumerate}
\end{rems}

The spectral assumption made in the above theorem are classically implied by a Kato condition on the Ricci curvature. Let $\left(P_t=e^{-t\Delta}\right)_{\{t\ge 0\}}$ be the heat semi-group and $H(t,x,y)$ the heat kernel, that is, for $f\in L^1(M)$ and $t>0$, one has 
$$\quad x\in M\colon\,(P_t f)(x)=\left(e^{-t\Delta}f\right)(x)=\int_M H(t,x,y)f(y)\dvol(y)  .$$
For a measurable  $V\colon M\to [0,\infty]$ and $T>0$, the Kato constant of $V$ is given by
$$\kappa_T(V):=\sup_{\ell \in\NN}\Vert  \int_0^T P_t(V\wedge \ell)\drm t\Vert_\infty,$$
(Note that the semigroup maps $L^\infty(M)$ to itself, such that the truncation procedure ensures that the quantity is well-defined.)
where $a\wedge b$ denotes the minimum of the numbers $a,b$. We say here that  a measurable $V\geq 0$ satisfies the Kato condition if there is a $T>0$ such that $$\kappa_T(V)<1.$$

The Kato condition is very powerful in proving, e.g., semiboundedness of the operator $\Delta-V$, and mapping properties of the corresponding semigroups. E.g., if the heat semigroup is ultracontractive, the perturbed semigroup will be, too. \\
\begin{rems}
There is another constant used in the literature also known as the Kato constant: for a measurable, non-negative $V\colon M \to[0,\infty]$, for some $L>0$ we let 
$$c_L(V):=\sup_{n\in\NN} \Vert (\Delta+L)^{-1}(V\wedge n)\Vert_\infty.$$ It is not hard to prove that $c_L$ and $\kappa_T$ are related by the formula \cite{Gueneysu-14}
$$(1-\euler^{-LT})c _L(V)\leq \kappa _T(V)\leq \euler^{LT}c _L(V),$$
meaning that the behavior of $\kappa_T(V)$ for $T\to 0$ controls the behavior of $c _L(V)$ for $L\to\infty$ and vice versa. 
\end{rems}

\begin{rems}
\begin{enumerate}
\item The Kato class was introduced in \cite{Kato-72} and further investigated in \cite{AizenmanSimon-82,Simon-82}, originally for the Laplacian in $\RR^n$. There, a function $V$ is in the Kato class if $c _T(V)\to 0$ as $T\to\infty$. Actually, the condition was phrased in terms of truncated heat kernels, see \cite{Voigt-86}.
\item The Kato condition has been extended to a certain class of measure perturbations of a general regular Dirichlet form in \cite{StollmannVoigt-96}. In this article, the Kato class consists of those measures $\mu$ for which $c _T(\mu)<1$ for large 
$T$. In particular, it was shown that in this context mapping properties of the semigroups carry over from the unperturbed to the perturbed. 
For instance, if for some $\upmu>0$, one has $\kappa_T(V)\le \left(1-e^{-\upmu T}\right)$, then
$$\Vert e^{-T(\Delta-V)}\Vert_{L^1\to L^1}\le e^{\upmu T}.$$ In particular, the bottom of the spectrum of the Schr\"odinger operator $\Delta-V$ is bounded from below by $-\upmu$.
\end{enumerate}
\end{rems}

\begin{coro}  If $(M^n,g)$ is a closed Riemannian manifold of dimension $n\ge 3$ such that for some $\uplambda>0$ and $T>0$ one has
$$\kappa_T\left( (\rho-\uplambda)_-\right)<\frac{1-e^{-(n-2)\uplambda T}}{n-2},$$ then the fundamental group of $M$ is finite.
\end{coro}

One can also obtain a Lichnerowicz-type estimate for the first non-zero eigenvalue. One of our results is the following:
\begin{propo}
Let $(M^n,g)$ be a closed manifold of dimension $n\ge 3$ such that for some $0\le k<1$, we have
 $$\kappa_T\left(\left(\rho -(n-1)\right)_-\right)\le \frac{n}{n-1}\left(1-e^{-T \frac{(n-1)^2}{n}(1-k^2)}\right).$$ Then, the first eigenvalue $\lambda_1(M)$ satisfies
$$\lambda_1(M)\geq  nk^2.$$
\end{propo}
In \cite{Carron-16}, it has been asked wether a control of the Ricci curvature in some Kato class yields some isoperimetric inequality. The second main result of this paper provided such a result using an idea of M.~Ledoux (\cite{Ledoux94}). The first result is a reverse Cheeger inequality in the spirit of P.~Buser (\cite{Buser82}). Recall that the Cheeger constant of $(M^n,g)$ is defined by:
$$\mathrm{h}(M,g):=\inf_\Omega \frac{\cA(\partial\Omega)}{\vol(\Omega)}, $$
where the infimum is taken over all $\Omega\subset M$ having smooth boundary $\partial\Omega$ and $\vol(\Omega)\leq \vol(M)/2$.
J.~Cheeger (\cite{Cheeger-70}) has shown that one always has $$\frac{\mathrm{h}^2(M,g)}{4}\le \lambda_1,$$ and P.~Buser proved that if the Ricci curvarture satisfies
$$\Ric_x\ge -(n-1)k^2,$$ then one has 
$$\lambda_1\le c(n)\left(k \mathrm{h}(M,g)+\mathrm{h}^2(M,g)\right).$$
\begin{propo}Let $(M^n,g)$ be a closed Riemannian manifold of dimension $n\geq 2$ such that for some $T>0$, we have
$$\kappa_T(\ricm)\le \frac 1{16n}\text.$$
There is a explicit constant $c_n$ such that  
$$\lambda_1(M)\le c_n\left( \frac{1}{\sqrt{T}}\mathrm h (M)+\mathrm h^2(M)\right).$$
\end{propo}
Using the eigenvalue estimates from \cite{Carron-16}, the estimate above yields lower bounds for the Cheeger constants under Kato-type assumptions on the negative part of the Ricci curvature. 
Note that assuming a lower bound $\Ric_x\ge -(n-1)k^2$ on the Ricci curvature yields $\kappa_T(\rho_-)\le T(n-1)k^2$, hence our result is a slight generalization of Buser's one. Furthermore, we obtain isoperimetric inequalities. Recall that for $p\in[1,\infty)$, the $p$-isoperimetric constant of $M$ is defined by
$$I_p(M):=\inf_\Omega \frac{\cA(\partial\Omega)}{\vol(\Omega)^{1-1/p}},$$
where the infimum is taken over all $\Omega\subset M$ having smooth boundary $\partial\Omega$ and $\vol(\Omega)\leq \vol(M)/2$.
\begin{theo} Let $(M^n,g)$ be a closed Riemannian manifold of dimension $n\geq 2$ and diameter at most $D>0$. 
Assume that for some $p>1$, there is an $\mathbf I>0$ such that
$$D^{2p-2}\kappa_T\left(\ricm^{p}\right)  \le \mathbf{I}^p,$$
and let
$\xi=\max\left\{ \frac{D}{\sqrt{T}}, \left(16n\mathbf{I}\right)^{\frac{p}{2p-2}}\right\}.$
Then 
we have
\begin{align*}
 D \vol(M)^{1/n}c_n^{1+\xi}I_n(M)\geq 1.
\end{align*}

\end{theo}

Our last result is that the assumptions made on the Kato constant of the Ricci curvature are implied under a control that depends on the volume of geodesic balls rather than the heat kernel.
\begin{theo} There is a positive constant $ \upeta_n$ such that when $(M^n,g)$ is a closed Riemannian manifold of diameter $D$ for which all $x\in M$ satisfy:
$$\int_0^{D} r \left(\frac{1}{\vol B(x,r)}\int_{B(x,r)} \ricm \right)dr\le \upepsilon_n$$
then $\kappa_{D^2}(\ricm)\le \frac{1}{16n}$ and the first Betti number of $M$ is smaller then $n$:
$$b_1(M)\le n.$$
\end{theo}
This paper is organized as follows. Sections~2.1 and 2.2 are devoted to present a Sobolev inequality under a spectral hypothesis on compact manifolds, which in turn implies a diameter estimate based on the ideas from \cite{BakrySobolev,Mondello2017}. We then prove a Lichnerowicz-type estimate under slightly stronger assumptions in Section~2.3. In Section~2.4 we prove a Bonnet-Myers theorem that is based on a spectral assumption of a certain Schr\"odinger operator. Section~3 presents all the Cheeger, Buser, and isoperimetric estimates, whereas in Section~4 we show how to bound the Kato constant under purely geometric assumptions.\\[1.5ex]
{\bf{Acknowledgement.}} G.C. was partially supported by the ANR grant CCEM-17-CE40-0034 and he wants to thank the Centre Henri Lebesgue ANR-11-LABX-0020-01 for creating an attractive mathematical environment. C.R. wants to thank the regional project D\'efiMaths and G.C. for their hospitality during his stay in Nantes, where parts of this work had been done. 
\section{A Sobolev inequality and diameter estimate under spectral hypotheses}
\subsection{Sobolev inequality}
The following is an elaboration from a result of D. Bakry \cite{BakrySobolev}  and revisited by I. Mondello \cite{Mondello2017}.
\begin{prop}\label{SobolevB} Let $(M^n,g)$ be a closed Riemannian manifold such that
for some $\delta\in \left(\frac{n+1}{n+4},1\right)$, the operator 
$$(1-\delta)\nabla^*\nabla+\ricci-(n-1)\Id$$ is non negative.\footnote{i.e. $\forall \alpha\in \cC^\infty(T^*M)\colon
(1-\delta)\int_M |\nabla \alpha|^2+\int_M\ricci(\alpha,\alpha)\ge (n-1)\int_M |\alpha|^2$}.
Then for $p_{n,\delta}=\frac{(1+\delta)n-1+\delta}{n-1-\delta}$ and $\gamma_{n,\delta}=\frac{(3+\delta)(n-\delta)}{(n-1-\delta)n(n-1)}$, we get the Sobolev inequality
$$\forall v\in \cC^\infty(M)\colon \|v\|_{L^{p_{n,\delta}}}^2\le \gamma_{n,\delta} \|dv\|^2_{L^2}+ \|v\|^2_{L^2}\ ;$$
where all norms are taken with respect to the probability measure
$d\mu=\frac{\dvol_g}{\vol_g(M)}.$
\end{prop}
 \proof We procced as in \cite{BakrySobolev}. For $p\in \left(2,\frac{2n}{n-2}\right)$ and $\upepsilon>0$ there is a positive function $f$ such that $\|f\|_{L^p}=1$ and that realizes the best constant in the inequality
 $$\forall v\in \cC^\infty(M)\colon \|v\|_{L^p}^2\le \gamma \|dv\|^2_{L^2}+ (1+\upepsilon)\|v\|^2_{L^2}.$$
 The Euler-Lagrange equation for the above problem implies that $f$ solves the equation
 $$\gamma\Delta f+(1+\upepsilon)f=f^{p-1}.$$
 For $\alpha\in \R$, we define $u$ by $f=u^\alpha$ and the computation of Bakry leads to
 \begin{equation}\label{equality1}
 \begin{split}
 (1+\upepsilon)\frac{p-2}{\gamma}\int_M |du|^2d\mu&=\int_M(\Delta u)^2d\mu\\
 &+(\alpha-1)(1+\alpha(p-2))\int_M \frac{|du|^4}{u^2}d\mu\\
 &-\alpha(p-1)\int_M \Delta u \frac{|du|^2}{u}d\mu.
 \end{split}\end{equation}
 With $A$ the traceless part of the Hessian of $u$: $$A=\nabla du+\frac{\Delta u}{n}g,$$
and the Bochner formula, we get
 \begin{equation*}
 \begin{split}
 \int_M(\Delta u)^2d\mu&=\int_M|\nabla d u|^2d\mu+\int_M\ricci(du,du)d\mu\\
&\ge (n-1)\int_M| d u|^2d\mu+\delta \int_M|\nabla d u|^2d\mu\\
&\ge (n-1)\int_M| d u|^2d\mu+\frac{\delta}{n} \int_M|\Delta u|^2d\mu+\delta \int_M|A|^2d\mu.
 \end{split}\end{equation*}
 Hence 
 $$\int_M(\Delta u)^2d\mu\ge \frac{n(n-1)}{n-\delta}\int_M| d u|^2d\mu+\frac{n\delta}{n-\delta} \int_M|A|^2d\mu.$$
 To estimate the last term of (\ref{equality1}), we integrate by parts and get
  \begin{equation*} \begin{split}
  \int_M \Delta u \frac{|du|^2}{u}d\mu=&\int_M \langle \nabla u, \nabla \frac{|du|^2}{u}\rangle d\mu\\
  &=\int_M  2\frac{\nabla du (du,du)}{u}d\mu-\int_M  \frac{|du|^4}{u^2} d\mu\\
  &=  -\frac{2}{n} \int_M \Delta u \frac{|du|^2}{u}d\mu+2\int_M\frac{A(du,du)}{u}d\mu-\int_M   \frac{|du|^4}{u^2}d\mu,
  \end{split}\end{equation*}
  hence one gets
  $$ \int_M \Delta u \frac{|du|^2}{u}d\mu=-\frac{n}{n+2} \int_M   \frac{|du|^4}{u^2}d\mu+\frac{2n}{n+2}\int_M\frac{A(du,du)}{u}d\mu.$$
  This implies
  \begin{equation}\label{inequality2}
 \begin{split}
 \left((1+\upepsilon)\frac{p-2}{\gamma}-\frac{n(n-1)}{n-\delta}\right)\int_M |du|^2d\mu=&\frac{n\delta}{n-\delta} \int_M|A|^2d\mu\\
 &-2n\alpha\frac{p-1}{n+2}\int_M\frac{A(du,du)}{u}d\mu\\
 &C(\alpha)  \int_M   \frac{|du|^4}{u^2}d\mu,\end{split}\end{equation}
where $C(\alpha)=(\alpha-1)(1+\alpha(p-2))+\frac{n}{n+2}\alpha(p-1).$
But since $\trace A=0$, we get
$$2\frac{A(du,du)}{u}=2\langle A, \frac{du\otimes du}{u} \rangle=2\langle A, \frac{du\otimes du}{u}-\frac{1}{n}\frac{|du|^2}{u}\Id \rangle$$
and 
$$2\left|\frac{A(du,du)}{u}\right|\le 2 |A|\sqrt{\frac{n-1}{n}} \frac{|du|^2}{u}\le \lambda |A|^2+\frac{n-1}{n\lambda}  \frac{|du|^4}{u^2}.$$
Chosing $\lambda$ such that 
$$\frac{n}{n-\delta}-n\alpha\frac{p-1}{n+2}\lambda=0,$$
we get
$$\left((1+\upepsilon)\frac{p-2}{\gamma}-\frac{n(n-1)}{n-\delta}\right)\int_M |du|^2d\mu\ge F(\alpha)\int_M   \frac{|du|^4}{u^2}d\mu,$$
where $F(\alpha)$ is the quadratic expression
$$F(\alpha)=A\alpha^2+B\alpha-1$$ with
$A=p-2\frac{(n-1)(p-1)^2(n-\delta)}{\delta(n+2)^2}$ and $B=2-2\frac{p-1}{n+2}$.
A little bit of arithmetic gives that 
$$\frac{B^2}{4}+A=\frac{n(p-1)}{n+2}\left(1+\frac{p-1}{n+2}\frac{-n+1+\delta}{\delta}\right).$$
This quantity is non positive if and only if 
$$0\le p-1\le \delta\frac{n+2}{n-1-\delta}.$$
Hence in that case, one can choose $\alpha$ such that $F(\alpha)\ge 0$ and we get
$$(1+\upepsilon)\frac{p-2}{\gamma}\ge \frac{n(n-1)}{n-\delta}.$$
It is then easy to conclude.
 \endproof
 \subsection{Diameter estimate} With the diameter estimates of Bakry-Ledoux \cite{BakryLedouxDiam}, we get 
 \begin{cor} Let $(M^n,g)$ be a closed Riemannian manifold satisfying the hypothesis of \pref{SobolevB}, then 
 $$\diam (M,g)\le \pi\, C(n,\delta)$$
 Where $C(n,\delta)=\frac{\sqrt{2p_{n,\delta}\upgamma_{n,\delta}}}{p_{n,\delta}-2}=\sqrt{\frac{1-(1-\delta)\frac{n-1}{n}}{1-(1-\delta)\frac{n+3}{4}}\left(1+\frac{1-\delta}{n-1}\right)}$.
 \end{cor}
 \begin{rems}
 \begin{enumerate}[i)]
 \item If $(M^n,g)$ is such  that for some $\delta>0$ with $1-\delta<\frac{4}{n+3}$, we have \begin{equation}\label{Positifkappa}
 (1-\delta)\nabla^*\nabla+\ricci\ge (n-1)k^2,\end{equation}
 then the metric $k^2g$ satisfies the hypothesis of the \pref{SobolevB} and hence the original metric $g$ satisfies the Sobolev type inequality
 $$\forall v\in \cC^\infty(M)\colon \|v\|_{L^{p_{n,\delta}}}^2\le \kappa \gamma_{n,\delta} \|dv\|^2_{L^2}+ \|v\|^2_{L^2},$$
 and hence the diameter estimate
 $$\diam M\le \frac{\pi}{k} C(n,\delta).$$
 \item 
If the bottom of the spectrum of the Schr\"odinger operator $(1-\delta)\Delta+\rho $ is strictly positive for some $\delta>0$ with $1-\delta<\frac{4}{n+3}$, i.e., 
 there is some $\mu>0$ with 
 $$\forall v\in \cC^\infty(M) \colon (1-\delta)\|dv\|^2_{L^2}+\int_M \rho  v^2\dvol_g\ge \mu\|v\|_{L^2},$$ then the condition (\ref{Positifkappa}) is satisfied with 
 $k=\sqrt{\mu/(n-1)}.$ 
 If this spectral condition is satisfied on $(M,g)$ then it is satisfied on any covering $\widehat M\rightarrow M$, but our result gives a Sobolev inequality and a diameter estimate only if $\widehat M$ is closed, i.e., for finite covering of $M$. We will see in the next section how such a spectral condition implies that the universal cover of $M$ has finite volume and hence that the fundamental group of $M$ is finite. \end{enumerate}
 \end{rems}
 One can give a Kato-type condition that implies the condition (\ref{Positifkappa}):
 Assume that $$\kappa_T((\rho -\uplambda)_-)\le (1-\delta)\left(1-e^{-\beta T}\right)$$
 where $0\le 1-\delta<\frac{4}{n+3}$. Then we get 
 $$\left\| e^{-T((1-\delta)\Delta-(\rho -\uplambda)_-}\right\|_{L^\infty\to L^\infty}\le e^{\frac{\beta}{1-\delta} T},$$
 and hence 
 $$(1-\delta)\Delta+\rho \ge -\frac{\beta}{1-\delta}+\uplambda.$$ If moreover $\beta<(1-\delta)\uplambda$, we get that the bottom of the spectrum of the Schr\"odinger operator $(1-\delta)\Delta+\rho $ is bounded from below by $-\frac{\beta}{1-\delta}+\uplambda$. This implies the following corollary:
 \begin{cor} Let $(M^n,g)$ be a closed manifold of dimension $n\ge 3$ such that for some $\upepsilon\in \left[0,\frac{4}{n+3}\right)$ and $k,\uplambda>0$ with $(n-1)k^2<\uplambda$, we have
 $$\kappa_T((\rho -\uplambda)_-)\le \upepsilon \left(1-e^{-\frac{T}{\upepsilon}(\uplambda-(n-1)k^2)}\right).$$ Then, the diameter of $M$ can be bounded by
 $$\diam(M)\le \frac{\pi}{k} C(n,1-\upepsilon). $$
 \end{cor}

\subsection{Lichnerowicz eigenvalue estimate}
It is easy to show that with the same arguments used in  the beginning of the proof of the \pref{SobolevB}, we get an eigenvalue estimate: 
\begin{lem}
Let $(M^n,g)$ be a closed manifold of dimension $n\geq 2$, such that for some $\delta\in (0,1]$ and $k>0$ the Schr\"odinger operator 
$(1-\delta)\Delta+\rho -(n-1)k^2$ is non negative. Then, we have 
$$\lambda_1(M)\geq \frac{n-1}{n-\delta}\, nk^2.$$ 
\end{lem}

We have found another method for proving a better estimate.
\begin{prop}Let $(M^n,g)$ be a closed manifold of dimension $n\geq 2$, such that for some  $k>0$ the Schr\"odinger operator 
$\frac{n}{n-1}\Delta+\rho -(n-1)k^2$ is non negative. Then, the first eigenvalue $\lambda_1(M)$ satisfies 
$$\lambda_1(M)\geq  nk^2.$$

\end{prop}

\proof By scaling, one can assume that $k=1$. 
 Let $\lambda>0$ be the first non-zero eigenvalue of the Laplacian with eigenfunction $\varphi$, i.e.,
 $$\Delta\varphi=\lambda\varphi.$$  For $\alpha>0$, we define $$h:=r^\alpha\varphi \ \in\cC^\infty((0,\infty)\times M)\text.$$ The manifold $K:=(0,\infty)\times M$, is equipped with the cone metric  $k=dr^2+r^2g$, such that the Laplacian $\Delta_h$ on $K$ has the representation $$\Delta_k=-\partial_r^2-\frac nr\partial_r+\frac 1{r^2}\Delta\text.$$
 We have
 \[
  \Delta_k h=\left(-\alpha(\alpha-1) -n\alpha+\lambda\right)r^{\alpha-2}\varphi\text.
 \]
 
Hence, we choose $\alpha>0$ so that $$\lambda:=\alpha(n-1+\alpha)$$  then $h$ will be a harmonic function on $(K,k)$. Note that $\alpha\ge 1$ is equivalent to $\lambda\ge n$.
 
 The Bochner formula on $K$ and the refined Kato inequality yield (see for instance \cite[Lemma 3.4]{BourCarron2})
$$\Delta_k\vert d_kh\vert^\frac{n-1}{n}\leq \frac{n-1}{n}\rho^K\vert d_kh\vert^\frac{n-1}{n},$$
where $d_k$ denotes the differential and $\rho^K$ is  the lowest eigenvalue of the Ricci tensor on $K$. Now, the norm on the cotangent bundle of $d_kh$ is
\[
 \vert d_kh\vert^2=r^{2\alpha-2}[\alpha^2\varphi^2+\vert d\varphi\vert^2],
\]
such that the second factor is independent of $r$. For $$\psi:=[\alpha^2\varphi^2+\vert d\varphi\vert^2]^\frac{n-1}{2n},$$
we have 
\[
 -(\alpha-1)\frac{n-1}{n}\left[(\alpha-1)\frac{n-1}{n}+n-1\right]\psi +\Delta\psi+\frac{n-1}{n}V\psi\leq 0,
\]
where $V(x)$ is the lowest eigenvalue of the Ricci curvature at $(1,x)$ in $(K,k)$.
A short calculation shows that $V(x)=\rho(x)-(n-1)$.
Hence, one gets
\begin{align}\label{proofLichnerowicz}
 \Delta\psi+\frac{n-1}{n}V\psi\leq (\alpha-1)\frac{n-1}{n}\left[(\alpha-1)\frac{n-1}{n}+n-1\right]\psi.
\end{align}
If we let
$$c(\alpha):=(\alpha-1)\frac{n-1}{n}\left[(\alpha-1)\frac{n-1}{n}+n-1\right],$$ then the above inequality \eqref{proofLichnerowicz} implies that the bottom of the spectrum of the operator $\Delta+\frac{n-1}{n}V$ is smaller than $c(\alpha)$. From our assumption, one gets $c(\alpha)\ge 0$, that is to say $\alpha\ge 1$.
\endproof

\begin{cor} Let $(M^n,g)$ be a closed manifold of dimension $n\ge 3$ such that for some $k,\uplambda>0$ with $\uplambda>(n-1)k^2$:
 $$\kappa_T((\rho -\uplambda)_-)\le \frac{n}{n-1}\left(1-e^{-T \frac{n-1}{n}(\uplambda-(n-1)k^2)}\right).$$ Then, $\lambda_1(M)$ satisfies
$$\lambda_1(M)\geq  nk^2.$$
\end{cor}

The Lichnerowicz-type estimate above also yields an analogue result for integral curvature bounds by estimating the norm of the heat semigroup as it was done in, e.g., \cite{RoseStollmann17,Rose-17, Gallot-88,Gallot-88_2}.

 \subsection{A Bonnet-Myers's theorem under a spectral hypothesis}
\begin{thm} If $(M^n,g)$ is a closed connected Riemannian manifold such that for some positive $\upepsilon\le 1/(n-2)$, the bottom of the spectrum of the Schr\"odinger operator
$$\upepsilon \Delta+\rho $$ is positive then $\uppi_1(M)$ is finite.
\end{thm}

\proof The main argument in the proof is in fact a modification of one part of the proof of the theorem in \cite{Carron-16}.
We are going to show that the universal cover $\pi\colon\widetilde M\rightarrow M$ has finite volume, hence the fundamental group is finite.

We first remark that if the bottom of the spectrum of the Schr\"odinger operator
$(n-2)\Delta+\rho$ is positive, then by continuity there is some positive $\upepsilon< 1/(n-2)$, such that the bottom of the spectrum of the Schr\"odinger operator
$\upepsilon\Delta+\rho$ is positive:
there is a positive constant $\mu$ such that $$\forall v\in \cC^\infty_0( M)\colon \upepsilon\int_{ M} |dv|^2+\int4_{M}\rho v^2\ge \mu\int_{M} v^2.$$
Let $\varphi\colon M\rightarrow \R$ be a positive eigenfunction of the Schr\"odinger operator $ \upepsilon\Delta+\rho$ associated to the bottom of the spectrum. We have $ \upepsilon\Delta\varphi+\rho\varphi=\lambda \varphi$ with $\lambda\ge \mu$. Hence on $\widetilde M$, one has for $\widetilde\varphi=\varphi\circ\pi$:
 $$ \upepsilon\Delta_{\widetilde M}\widetilde\varphi+\rho\widetilde\varphi=\lambda \widetilde\varphi.$$   By a principle due to W.F.~Moss and J.~Piepenbrink, and D.~Fisher-Colbrie and R. Schoen, (\cite{MossPiep-78,FischerSchoen-80} or \cite[lemma 3.10]{Pigola-Book-08}), that is in fact based on an old result of J.~Barta (\cite{Barta-37}), we know that  $(\widetilde M,\tilde g=\pi^*g)$ also satisfies
$$\forall v\in \cC^\infty_0(\widetilde M)\colon \upepsilon\int_{\widetilde M} |dv|^2+\int_{\widetilde M}\rho v^2\ge \mu\int_{\widetilde M} v^2.$$
With the Bochner formula and the Kato inequality, we get that the spectrum of the Hodge-deRham Laplacian on one-forms is positive:
$$\spec \Delta_1\subset [\mu,+\infty).$$
Indeed, if $\alpha\in \cC_0^\infty(T^*\widetilde M)$ then
\begin{equation*}\begin{split}
\langle \Delta_1\alpha,\alpha\rangle&=\int_{\widetilde M} |\nabla\alpha|^2+\ricci(\alpha,\alpha)\\
&\ge \int_{\widetilde M} |d|\alpha||^2+\rho |\alpha|^2\\
&\ge \mu \int_{\widetilde M} |\alpha|^2.
\end{split}\end{equation*}
Because one always has $\spec \Delta\subset \{0\}\cup \spec \Delta_1$, one deduces that the spectrum of the Laplacian $\Delta$ on functions has a spectral gap:
$$\spec \Delta\subset \{0\}\cup [\mu,+\infty).$$

Hence either $0$ is an eigenvalue for the Laplacian on $\widetilde M$ and this means that the volume of $\tilde M$ is finite\footnote{Indeed, if $f$ is an eigenfunction of the Laplace operator associated to the eigenvalue $0$, then $\Delta f=0$ and $df\in L^2$ hence by definition of the domain of the Laplace operator, one can integrate by parts and gets
$$\int_{\widetilde M}|df|^2=\int_{\widetilde M}f \Delta f =0.$$ Hence $f$ is constant. } or the bottom of the spectrum of the Laplacian is positive. It is classical that the first case implies that $\uppi_1(M)$ is finite. We are going to show that the second case can not occur. 

We are arguing by contradiction and we suppose that $\lambda_0(\widetilde M)>0$. 
$\widetilde M$ being a covering of a compact Riemannian manifold, it satisfies the local Sobolev inequality :
$$\forall v\in \cC^\infty_0(\widetilde M)\colon \|v\|_{L^{\frac{2n}{n-2}}}^2\le C\left( \|dv\|^2_{L^2}+ \|dv\|^2_{L^2}\right).$$
Using the spectral gap, we get the Euclidean Sobolev inequality:
$$\forall v\in \cC^\infty_0(\widetilde M)\colon \|v\|_{L^{\frac{2n}{n-2}}}^2\le C\left(1+\lambda_0(\widetilde M)^{-1}\right) \|dv\|^2_{L^2}$$
In particular, $(\widetilde M,\tilde g)$ is non-parabolic and there is a unique positive minimal Green kernel $G(x,y)$.
Let $o\in \widetilde M$ be a fixed point and define $b\colon M\rightarrow \R_+$ by
$$ b^{2-n}=c_n G(o,x),$$ where $c_n$ is chosen so that 
$b(x)\simeq d(o,x)$ as $x\to o$. The Sobolev inequality implies that there is a positive constant $c$ such that 
$$c\ d_{\tilde g}(o,x)\le b(x).$$
In particular, $\lim_{x\to \infty} b(x)=+\infty$.
According to T.H.~Colding (\cite{ColdingActa}), for any $\upalpha\ge (n-2)/(n-1)$, we have
$$\Delta \frac{|db|^\upalpha}{b^{n-2}}+\upalpha\rho   \frac{|db|^\upalpha}{b^{n-2}}\le 0\ \mathrm{weakly\ on\ } \widetilde M\setminus \{o\}.$$
Our hypothesis implies that for any $\upalpha\in [0,1/\upepsilon]$, the Schr\"odinger operator $\Delta+\upalpha \rho $ has a spectral gap, hence it has a unique minimal positive Green kernel $G_\upalpha$. We let $g_\upalpha(x)=G_\upalpha(o,x)/c_n.$ It is then classical (see \cite[subsection 2.6.2]{Carron-16}) that 
\begin{equation}\label{GreenH1}\int_{\widetilde M\setminus B(o,1)} |dg_\upalpha|^2<\infty.\end{equation}
Then the spectral gap $$\forall \varphi\in \cC^\infty(\widetilde M)\colon \int_{\widetilde M} \varphi^2\le \int_{\widetilde M} |d\varphi|^2+\alpha\rho\varphi^2$$ and the proof of the above estimate (\ref{GreenH1}) implies that 
$$\int_{\widetilde M\setminus B(o,1)} G_\upalpha^2<\infty.$$

Note that with the local Sobolev inequality, the classical De Giorgi-Nash-Moser iteration scheme (see \cite[theorem 8.17]{gilbarg2015elliptic} or  \cite[Lemma 1.5]{Mondello2017}) implies that a weak non-negative $W^{1,2}_{\mathrm{loc}}$ solution of 
$$\Delta\upvarphi \le \Lambda \upvarphi$$ is in fact locally bounded: for any $p\ge 1$, there is a constant $C$ depending only on $r,\Lambda$ and the local Sobolev constant so that 
$$\sup_{y\in B(x,r)}\upvarphi(y)\le C\left(\int_{B(x,2r)} \upvarphi^p\right)^{\frac1p}.$$

$(\widetilde M,\tilde g)$ being a covering of a compact Riemannian manifold, the Ricci curvature is bounded from below: $\ricci_{\tilde g}\ge -(n-1)k^2$ so that $-\rho\le (n-1)k^2$. Hence  we get the following elliptic estimate: 
for any $\upalpha$ and $p\ge 1$ there is a constant $C$ such that if a non negative function $\upvarphi\in W^{1,2}_{loc}$ satisfies 
$$\Delta\upvarphi +\upalpha\rho   \upvarphi\le 0\ \mathrm{weakly\ on\ } B(x,2)$$
then 
$$\sup_{y\in B(x,1)}\upvarphi(y)\le C\left(\int_{B(x,2)} \upvarphi^p\right)^{\frac1p}.$$
In particular, we get (when $\upalpha\in [0,1/\upepsilon]$) that $\lim_{x\to \infty} g_\upalpha(x)=0$,
and we have the equivalence 
$$\frac{|db|^\upalpha}{b^{n-2}}\le g_\upalpha\Longleftrightarrow \lim_{x\to \infty} \frac{|db|^\upalpha}{b^{n-2}}=0.$$
With the estimate (\ref{GreenH1}), one gets $$\int_{\widetilde M\setminus  B(o,1)} |d_xG(o,x)|^2\dvol_{\tilde g}(x)<\infty,$$ we deduce that for $\upalpha_0:=\frac{n-2}{n-1}$, we have
$$\lim_{x\to \infty} \frac{|db|^{\upalpha_0}}{b^{n-2}}=0,$$ and hence $$\frac{|db|^{\upalpha_0}}{b^{n-2}}\le g_{\upalpha_0}.$$

Our main tool is the universal Hardy inequality (\cite{CarronHardy}):
$\forall  \psi\in \cC_0^{\infty}(\widetilde M)$:
$$\frac{(n-2)^{2}}{4}\int_{\widetilde M} \frac{|db|^{2}}{b^{2}}\psi^{2}\dvol_{\tilde g}\le \int_{\widetilde M}|d\psi|^{2}\dvol_{\tilde g}.$$
The \sch operator $\upepsilon \Delta+\rho $ is positive hence  for any $\upalpha\in [0,n-2]$, we have the following Hardy type inequality:
$\forall  \psi\in \cC_0^{\infty}(\widetilde M)$:
$$c\int_{ \widetilde M} \frac{|db|^{2}}{b^{2}}\psi^{2}\dvol_{\tilde g}\le \int_{ \widetilde M}\left[|d\psi|^{2}+\upalpha\rho  \psi^2\right]\dvol_{\tilde g}.$$
When $\upalpha\in [\upalpha_{0},n-2]$, using the function $\psi=\xi \frac{|db|^{\upalpha}}{b^{n-2}}$ where $\xi$ is a Lipschitz function with compact support in $\widetilde M\setminus \{o\}$ one gets
\begin{equation*}
c\int_{\widetilde M} \frac{|db|^{2+2\upalpha}}{b^{2(n-1)}}\xi^{2}\dvol_{\tilde g}\le \int_{\widetilde M}|d\xi|^{2}\frac{|db|^{2p}}{b^{2(n-2)}}\dv_{g}.
\end{equation*}

Assume that for some $\upalpha \in [\upalpha_{0},n-2]$, we have
$$\frac{|db|^{\upalpha}(x)}{b^{n-2}(x)}\le g_{\upalpha}(o,x).$$
Then we get 
\begin{equation*}\forall \xi\in \cC^\infty_0(\widetilde M\setminus\{o\})\colon
c\int_{\widetilde M} \frac{|db|^{2+2\upalpha}}{b^{2(n-1)}}\xi^{2}\dvol_{\tilde g}\le \int_{\widetilde M}|d\xi|^{2} g_\upalpha^2\dvol_{\tilde g}.
\end{equation*}
According to \cite[Proposition 2.21]{Carron-16}, if $\xi$ is a smooth function with support in $\widetilde M\setminus B(o,1/2)$ that is identitically $1$ on $\widetilde M\setminus B(o,1)$, then there is a sequence of smooth functions $(\xi_\ell)_\ell$ with compact support in $\widetilde M\setminus\{o\}$ such that 
$$\lim_\ell \xi_\ell=\xi\ \mathrm{uniformly\ on\ compact\ subset\ of\ }\widetilde M.$$
and 
$$\lim_\ell  \int_{\widetilde M}|d\xi_\ell|^{2} g_\upalpha^2\dv_{g}= \int_{\widetilde M}|d\xi|^{2} g_\upalpha^2\dv_{g}.$$
In particular, we get 
$$\int_{\widetilde M\setminus B(o,1)} \frac{|db|^{2+2\alpha}}{b^{2(n-1)}}\dv_{g}<\infty$$
If $\upalpha'=(1+\upalpha)\frac{n-2}{n-1}$ one gets
$$\int_{M\setminus B(o,1)} \left(\frac{|db|^{\upalpha'}}{b^{n-2}}\right)^{2\frac{n-1}{n-2}}\dv_{g}<\infty,$$
hence, 
$$\lim_{x\to \infty} \frac{|db|^{\upalpha'}(x)}{b^{n-2}(x)}=0$$ and 
$$\frac{|db|^{\upalpha'}(x)}{b^{n-2}(x)}\le g_{\upalpha'}.$$
For $\upalpha_{k}=(n-2)-\left(\frac{n-2}{n-1}\right)^{k}\left(n-2-\upalpha_{0}\right)=\left(1+\upalpha_{k-1}\right)\frac{n-2}{n-1},$ our argumentation yields that for all $k\in \N$:
$$\frac{|db|^{\upalpha_{k}}(x)}{b^{n-2}(x)}\le g_{ \upalpha_{k}}(o,x).$$
Hence letting $k\to\infty$, we obtain the following estimate on the log derivative of the Green kernel:
$$\frac{|db|^{n-2}(x)}{b^{n-2}(x)}\le g_{ n-2}(x).$$
And in particular $$\lim_{x\to \infty} \frac{|db|(x)}{b(x)}=0.$$
We are going to show that this implies that $\lambda_0(\widetilde M)=0$, hence a contradiction.
The Green formula implies that if $t$ is a regular value of $b$, then\footnote{with $\sigma_{n-1}=\vol \bS^{n-1}$.}
$$\int_{\{b=t\}} |db|=\sigma_{n-1}t^{n-1},$$ 
so that\footnote{With $\omega_n=\vol \bB^n$.}
$$\int_{\{b\le R\}} |db|^2=\omega_n R^n$$ and
$$\int_{\{R\le 4b\le 3 R\}} |db|^2=\omega_n (3^n-1)R^n.$$
Let $\Omega_R:=\{b\le R\}$ and $\upvarphi_R=(R-b)_+$. We have that $\upvarphi_R$ has support in $\Omega_R$ and
$$\int_{\Omega_R} |d\upvarphi_R|^2=\omega_n R^n=\frac{1}{3^n-1} \int_{\{R\le 4b\le 3 R\}} |db|^2.$$
But 
$$\int_{\Omega_R} \upvarphi_R^2\ge  \int_{\{R\le 4b\le 3 R\}}  \upvarphi_R^2\ge \frac{1}{9}\int_{\{R\le 4b\le 3 R\}}  b^2.$$
We have shown that $$\lim_{x\to \infty} \frac{|db|(x)}{b(x)}=0,$$ hence 
$$\lim_{R\to\infty} \frac{\int_{\Omega_R} |d\upvarphi_R|^2}{\int_{\Omega_R} |\upvarphi_R|^2}=0$$ hence
$\lambda_0(\widetilde M)=0$.
\endproof

\section{Buser inequality and isoperimetric estimates}

We are going to use an idea from M. Ledoux \cite{Ledoux94} in order to obtain isoperimetric inequalities. 
\subsection{Buser inequality}
To apply Ledoux's technique, we recall the gradient estimates for positive solutions of the heat equation from \cite{Rose:2016aa,Carron-16}:
\begin{prop}\label{LiYau}
Let $(M^n,g)$ be a closed Riemannian manifold of dimension $n\geq 2$ such that for some $T>0$ :
$$\kappa_T(\ricm)\le \frac 1{16n}\text.$$
If $u\colon [0,T]\times M\rightarrow \R_+$ be a positive solution of the heat equation
\begin{align}\label{heatequation}
\partial_tu=-\Delta u\text,
\end{align}

then, we have
\begin{align}\label{gradientestimate}
e^{-2}\frac{\vert \nabla u\vert^2}{u^2}-\frac{\partial_t u}{u}\leq \frac{e^2n}{ 2t }\text, \quad \forall \, t\in(0,T)\text.
\end{align}
\end{prop}

From Proposition \ref{LiYau}, we can deduce
\begin{cor} Under the assumptions of \pref{LiYau}, there is an explicit constant $c_n$ depending only on the dimension, such that for any $t\in (0,T]$:
$$\|\Delta P_t\|_{L^\infty\to L^\infty}\le \frac{c_n}{t},$$
and
$$\|\nabla P_t\|_{L^\infty\to L^\infty}\le \frac{c_n}{\sqrt{t}}.$$
\end{cor}
\proof
Let $f\in L^1(M)$ be a non negative initial value for the heat equation with $u=P_tf$. 
Note that 
\begin{align*}(\Delta u)_+=\left(-\partial_tu\right)_+\le  \frac{e^2n}{ 2t } u.
\end{align*}
By Stokes theorem, we have $\int\Delta u=0$. Furthermore, we know $P_t1=1$ and $\int P_tf=\int f$. Hence,

\begin{align*}
 \Vert\Delta P_tf\Vert_{L^1}=2\Vert(\Delta P_tf)_+\Vert_{L^1}\leq  \frac{e^2n}{ t }\Vert u\Vert_{L^1}= \frac{e^2n}{ t }\Vert f\Vert_{L^1}.
\end{align*}
We easily get 
$$\|\Delta P_t\|_{L^1\to L^1}\le \frac{e^2n}{t}.$$
By duality and self-adjointness of the operator $\Delta P_t$ we obtain the first assertion.

Let $f\in L^\infty(M)$ be a non negative initial value for the heat equation and let $u=P_tf$. We have
\begin{align*}\vert \nabla u\vert^2\le e^{2}\left(\frac{e^2n}{ 2t } u^2+u\left| \Delta u\right| \right),
\end{align*}
and, therefore,
\begin{align*}\Vert \nabla u\Vert^2_{L^\infty}&\le e^{2}\left(\frac{e^2n}{ 2t }+\frac{e^2n}{t}\right)\Vert  u\Vert^2_{L^\infty}\\
&\le\frac{3e^4n}{ 2t }  \Vert  f\Vert^2_{L^\infty}.\end{align*}
In general, for $f\in L^\infty(M)$, we get
\begin{align*}\Vert \nabla P_t f\Vert_{L^\infty}\le \sqrt{\frac{14n}{t}}\Vert  f\Vert_{L^\infty}.\end{align*}
\endproof
The above proposition implies the following crucial estimate. 
\begin{prop}\label{pseudoPoinc} Under the assumptions of \pref{LiYau}, there is a constant\footnote{One can choose $c_n=6\sqrt{n}$.} $c_n$, such that for  any $t\in (0,T]$, we have
\begin{align*}\forall f\in W^{1,1}(M)\colon \Vert f-P_tf\Vert_{L^1}\le c_n \sqrt{t}\ \Vert df\Vert_{L^1}.
\end{align*}
\end{prop}
\proof
We have 
$$ f-P_tf=\int_0^t \Delta e^{-s\Delta}fds.$$
For any $g\in L^\infty$, the above yields
\begin{align*}
 \int_M g\Delta e^{-s\Delta}f\dvol&=\int_M  g\Delta P_sf\dvol\\
 &=\int_M \langle\nabla P_t g,\nabla f\rangle \dvol\\
  &\leq \sqrt{\frac{14n}{t}}\Vert\nabla f\Vert_{L^1}\Vert g\Vert_{L^\infty}.
\end{align*}
Hence, we get
\begin{align*}\Vert\Delta e^{-s\Delta}f\Vert_{L^1}\le \sqrt{\frac{14n}{t}}\Vert\nabla f\Vert_{L^1},\end{align*}
what implies
\begin{align*}\Vert f-P_tf\Vert_{L^1}\le c_n \sqrt{t} \Vert df\Vert_{L^1},\end{align*}

with $c_n=6\sqrt{n}.$
\endproof
Now, we can obtain the following estimate relating the Cheeger constant and the first non-zero eigenvalue of the Laplacian on $M$. 
\begin{thm}\label{katobuser}Let $(M^n,g)$ be a closed Riemannian manifold of dimension $n\geq 2$ such that for some $T>0$, we have
$$\kappa_T(\ricm)\le \frac 1{16n}\text.$$
There is a explicit constant $c_n$ such that  if $\lambda_1$ is the first non-zero eigenvalue of the Laplacian on $M$
satisfies
$$\lambda_1(M)\le c_n\left( \frac{1}{\sqrt{T}}\mathrm h(M)+\mathrm h^2(M)\right).$$
\end{thm}
\proof We apply the above inequality to the characteristic function $\chi_{\Omega}$ of $\Omega$.
From the fact that $P_t(\chi_\Omega)\le 1$ we get
\begin{align*}
\int_M \left|\chi_{\Omega}-P_t(\chi_{\Omega})\right| \dvol&=\int_\Omega \left|\chi_{\Omega}-P_t(\chi_{\Omega})\right| \dvol+\int_{M\setminus\Omega} P_t(\chi_{\Omega})\dvol\\
&\int_\Omega \left(\chi_{\Omega}-P_t(\chi_{\Omega})\right) \dvol+\int_{M\setminus\Omega} P_t(\chi_{\Omega})\dvol.
\end{align*}
Since
$\int_\Omega P_t(\chi_{\Omega})\dvol=\Vert P_{t/2}\chi_\Omega\Vert_{L^2}^2$ and
$$\int_{M\setminus\Omega} P_t(\chi_{\Omega})\dvol=\langle P_t(\chi_{M\setminus \Omega}),\chi_{\Omega}\rangle=\langle 1-P_t(\chi_{ \Omega}),\chi_{\Omega}\rangle=\vol \Omega-\Vert P_{t/2}\chi_\Omega\Vert_{L^2}^2,$$
we have
\begin{align*}
\int_M \left|\chi_{\Omega}-P_t(\chi_{\Omega})\right| \dvol=2\left(\vol \Omega-\Vert P_{t/2}\chi_\Omega\Vert_{L^2}^2\right).\end{align*}
We can decompose $\chi_\Omega=\frac{\vol\Omega}{\vol M}+f$ where $\int_M f\dvol=0$ and get
$$\Vert P_{t/2}\chi_\Omega\Vert_{L^2}^2\le \frac{(\vol\Omega)^2}{\vol M}+e^{-\lambda_1 t} \left(\vol\Omega- \frac{(\vol\Omega)^2}{\vol M}\right).$$
This yields
\begin{align*}
2\vol \Omega  \left(1- \frac{\vol\Omega}{\vol M}\right)\left(1- e^{-\lambda_1 t}\right)\le c_n\sqrt{t}\vol \partial \Omega,
\end{align*}
and we get for any $t\in [0,T]$ that
$$\left(1- e^{-\lambda_1 t}\right)\le c_n\sqrt{t}\ \mathrm h(M).$$

If $T\lambda_1\ge 1$, then we choose $t=1/\lambda_1$ and get
$$\sqrt{\lambda_1} \le 3c_n \mathrm h(M).$$
In the case where $T\lambda_1\le 1$, we choose $t=T$ and from the inequality
$2(1-e^{-x})\geq x$ for $x\in (0,1)$, we get
$$\lambda_1 \sqrt{T}\le 2c_n \mathrm h(M).$$
It is now easy to get the announced result.
\endproof
\subsection{An estimate of the Cheeger constant}
The following corollary is an easy consequence of Theorem~\ref{katobuser}.
\begin{cor}Let $(M^n,g)$ be a closed Riemannian manifold of dimension $n\geq 2$ and diameter at most $D>0$.  There is a constant $c_n\ge 1$ such that if   $T$ is the first positive time such that 
$$\kappa_T(\ricm)= \frac 1{16n},$$ and if we define $\xi:=\xi(M,g)=D/\sqrt{T}$,
then 
$$\frac{c_n^{-1-\xi}}{D}\le \mathrm h(M).$$
\end{cor}
\proof As a matter of fact, with the eigenvalue estimate given in \cite[Theorem3.6]{Carron-16}, we know that 
$$\frac{c_n^{-1-\xi}}{D^2}\le \lambda_1.$$ Therefore, we have
$$\frac{c_n^{-1-\xi}}{D^2}\le c_n\left( \frac{\xi }{D}\mathrm h(M)+\mathrm h^2(M)\right).$$
Hence either
$\frac{1}{D}\le h(M)$ or
$$\frac{c_n^{-1-\xi}}{D^2}\le c_n\frac{1+\xi }{D}\mathrm h(M)\le c_n \frac{e^\xi}{D} \mathrm h(M). $$
\endproof

\subsection{Isoperimetric inequalities}
From the results above, we can now obtain isoperimetric inequalities.
\begin{thm}\label{isoperimetryprofile} Let $(M^n,g)$ be a closed Riemannian manifold of dimension $n\geq 2$ and diameter at most $D>0$. 
Assume one of the following:
\begin{enumerate}[i)]
\item $\kappa_T(\ricm)\le \frac{1}{16n},$ 
\item or that for some $p>1$, we have for some $\mathbf I>0$ that 
$D^{2p-2}\kappa_T\left(\ricm^{p}\right)  \le \mathbf{I}^p.$
\end{enumerate}
In the first case let $\xi:=\xi(M,g)=D/\sqrt{T}$ and $\nu=ne^{8\sqrt{n\kappa_T(\ricm)}}$ and in the second case let
$\xi=\max\left\{ \frac{D}{\sqrt{T}}, \left(16n\mathbf{I}\right)^{\frac{p}{2p-2}}\right\},$ and $\nu=n$.

Then 
we have
\begin{align*}
1\leq c_n^{1+\xi}\vol(M)^\frac{1}{\nu}D\, I_\nu(M).
\end{align*}
\end{thm}
\proof
From \cite[Proposition 3.7 and 3.14]{Carron-16}, we get for any $t\in \left[0,D^2\right]$:
\begin{align*}\Vert P_t\Vert_{L^1\to L^\infty}\le \frac{c_n^{1+\xi}}{\vol M}\ \frac{D^\nu}{t^{\frac \nu 2}}.
\end{align*}
Let $\Omega\subset M$ with smooth boundary and $\vol \Omega\le \frac12 \vol M$.
Then, we have
\begin{align*}
\vol \Omega\le \vol\left\{x\in M, |\chi_\Omega-P_t(\chi_\Omega)|>1/2\right\}+ \vol\left\{x\in M, P_t(\chi_\Omega)>1/2\right\}.
\end{align*}
Furthermore,
\begin{align*}\Vert P_t(\chi_\Omega)\Vert_{L^\infty}\le \frac{c_n^{1+\xi}\vol \Omega}{\vol M}\ \frac{D^\nu}{t^{\frac \nu 2}}.
\end{align*}
Hence if 
$$\frac{c_n^{1+\xi}\vol \Omega}{\vol M}\le \frac{1}{4}$$ one can choose $t\in [0,D^2]$ such that 
$$\frac{c_n^{1+\xi}\vol \Omega}{\vol M}\ \frac{D^\nu}{t^{\frac \nu 2}}=\frac 14,$$
such that

\begin{align*}
\vol \Omega\le 2\Vert \chi_\Omega-P_t(\chi_\Omega)\Vert_{L^1}\le c_n \sqrt{t} \vol \partial \Omega.\end{align*}
A little bit of arithemetics then implies the result when $\frac{c_n^{1+\xi}\vol \Omega}{\vol M}\le \frac{1}{4}$. When
$\frac{c_n^{1+\xi}\vol \Omega}{\vol M}\ge \frac{1}{4}$, we use the estimate on the Cheeger constant to conclude.
\endproof
\section{An estimate of the Kato constant }
Let $(M^n,g)$ be a complete Riemannian manifold. We assume that it satisfies the following volume doubling condition and heat kernel estimates at scale $R$:
there are positive constants  $\uptheta,\nu,C$  such that 
\begin{equation}\label{doubling}\tag{D${}_R$} \forall x\in M, 0\le r\le \rho\le R\colon \frac{\vol B(x, \rho)}{\vol B(x,r)}\le \uptheta\left(\frac{\rho}{r}\right)^{\nu}\end{equation}
\begin{equation}\label{heatkernel}\tag{UE${}_R$} \forall x,y\in M, t\in (0,R^2]\colon H(t,x,y)\le \frac{C}{\vol B(x,\sqrt{t})}e^{-\frac{d^2(x,y)}{5t}}.\end{equation}

Recall that if $q\colon M\rightarrow \R_+$, is locally integrable then its parabolic Kato constant at times $T$ is given by
$$\kappa_T(q)=\sup_{x\in M} \int_0^T\int_M H(t,x,y)q(y)\dvol_g(y)dt$$
We know that if $\kappa_T(q)\le 1-e^{-\upbeta T}$ for some $\upbeta>0$, then 
$$\forall t>0\colon\ \left\|e^{-t(\Delta-q)} \right\|_{L^\infty\to L^\infty}\le e^{\upbeta t},$$ and in particular the Schr\"odinger operator $\Delta-q$ is bounded from below by $-\upbeta$:
$$\forall v\in \cC^\infty_0(M)\colon \int_M \left[ |dv|^2-qv^2\right]\dvol_g\ge -\upbeta \int_Mv^2\dvol_g.$$
\begin{prop} There is a constant $\uplambda$ depending only on $\theta,\nu$ and $C$ such that if \begin{equation}\label{condition}
\lim_{r\to +\infty}e^{-\frac{r^2}{R2}}\int_{B(x,r)} q=0\end{equation} 
then
$$\kappa_{R^2}(q)\le \uplambda\sup_x\int_0^\infty r e^{-\frac{r^2}{6R^2}}\left(\fint_{B(x,r)} q\right)dr.$$
\end{prop}
\proof We introduce the convenient notations :
$$V(x,r)=\vol B(x,r)\ \mathrm{and}\ Q(x,r)=\int_{B(x,r)} q$$ so that 
$$\fint_{B(x,r)} q=\frac{Q(x,r)}{V(x,r)}.$$
Notice that $r\mapsto Q(x,r)$ is a non-decreasing function so that we can consider the Riemann-Stieljes measure
$dQ(x,r)$.
Using the estimate on the heat kernel we have 
$$\int_0^{R^2}\int_M H(t,x,y)q(y)\dvol_g(y)dt\le C\int_{[0,R^2]\times [0,\infty)} \frac{e^{-\frac{r^2}{5t}}}{V(x,\sqrt{t})} dQ(x,r)dt.$$
The condition (\ref{condition}) implies  that we can integrate by parts:
$$\int_{[0,R^2]\times [0,\infty)} \frac{e^{-\frac{r^2}{5t}}}{V(x,\sqrt{t})} dQ(x,r)dt=\int_{[0,R^2]\times [0,\infty)} \frac{2r e^{-\frac{r^2}{5t}}}{5t V(x,\sqrt{t})} Q(x,r)dtdr.$$
We are going to split this integral in 3 parts corresponding to the domains $\{r\ge R\}$ $\{\sqrt{t} \le r\le R\}$ and $\{ r\le \sqrt{t} \le R\}$.

Concerning the first part, when $r\ge R$, we get
$$\int_0^{R^2}\frac{e^{-\frac{r^2}{5t}}}{t V(x,\sqrt{t})} dt\le \uptheta\frac{R^\nu}{V(x,R)}\int_0^{R^2} \frac{e^{-\frac{r^2}{5t}}}{t^{\frac{\nu}{2}+1}} dt.$$
Using $\xi=r^2/(5t)$, we obtain
$$\int_0^{R^2} \frac{e^{-\frac{r^2}{5t}}}{t^{\frac{\nu}{2}+1}} dt\le \frac{5^\nu}{r^\nu}\int_{\frac{r^2}{5R^2}}^\infty e^{-\xi}\xi^{\frac{\nu}{2}-1}d\xi.$$
Now there is a constant $c(\nu)$ such that when $A\ge 1/5$:
$$\int_{A}^\infty e^{-\xi}\xi^{\frac{\nu}{2}-1}d\xi\le c(\nu) e^{-A} A^{\frac{\nu}{2}-1}$$ and we get
$$\int_0^{R^2}\frac{e^{-\frac{r^2}{5t}}}{t V(x,\sqrt{t})} dt\le \uptheta c(\nu)5^\nu \frac{e^{-\frac{r^2}{5R^2}}}{V(x,R)}.$$
According to \cite[Lemma 3.10]{Carron-16}, the local doubling condition implies that for $r\ge R$:
$$V(x,r)\le (\uptheta 2^\nu)^{50+50\frac{r}{R}}V(x,R)$$
Hence we get that for some constant $\uplambda_1$ depending only on $\uptheta,\nu, C$:
$$\int_0^{R^2}\int_{M\setminus B(x,R)} H(t,x,y)q(y)\dvol_g(y)dt\le \uplambda_1 \int_{R}^\infty e^{-\frac{r^2}{5R^2}}\frac{r}{V(x,r)}Q(x,r)dr.$$

For the second part, we easily get by the same argumentation:
$$\int_0^{r^2}\frac{e^{-\frac{r^2}{5t}}}{t V(x,\sqrt{t})} dt\le \uptheta\frac{r^\nu}{V(x,r)}\int_0^{r^2} \frac{e^{-\frac{r^2}{5t}}}{t^{\frac{\nu}{2}+1}} dt\le c(\nu) \uptheta\frac{1}{V(x,r)}.$$
For the remaining part, we have to estimate:
\begin{equation*}\begin{split}
\int_0^{R}\left(\int_{r^2}^{R^2} \frac{r e^{-\frac{r^2}{5t}}}{t V(x,\sqrt{t})}dt\right) Q(x,r)drdt&=\int_0^{R^2}\frac{1}{V(x,\sqrt{t})}\left(\int_0^{\sqrt{t}}\frac{r e^{-\frac{r^2}{5t}}}{t }Q(x,r)dr\right)dt\\
&\le \int_0^{R^2}\frac{Q(x,\sqrt{t})}{V(x,\sqrt{t})}\left(\int_0^{\sqrt{t}}\frac{r e^{-\frac{r^2}{5t}}}{t }dr\right)dt\\
&\le \int_0^{R^2}\frac{Q(x,\sqrt{t})}{V(x,\sqrt{t})} \left[-\frac{5}{2}e^{-\frac{r^2}{5t}}\right]_0^{\sqrt{t}} dt\\
&\le  \frac{5}{2}\int_0^{R^2}\frac{Q(x,\sqrt{t})}{V(x,\sqrt{t})}dt= 5\int_0^{R}\frac{Q(x,r)}{V(x,r)}rdr
\end{split}\end{equation*}

\endproof
\begin{cor}\label{estikatoclosed} If $(M^n,g)$ is a closed Riemannian manifold with diameter $D$ and that satisfies the doubling (\ref{doubling}) and heat kernel estimate (\ref{heatkernel}) at scale $R$ with $R\le D$, then 
$$\kappa_{R^2}(q)\le \uplambda\sup_x\int_0^D r e^{-\frac{r^2}{7R^2}}\left(\fint_{B(x,r)} q\right)dr.$$
\end{cor}
\proof We know that there is a positive constant $C$ such that if $\gamma=e^{C\frac{D}{R}}$ then
$$\forall x\in M\colon \vol_g(M)\le \gamma \vol B(x,D/4).$$
Hence there is a finite set $A\subset M$ with $\#A\le \gamma$ such that 
$$M=\cup_{a\in A} B(a,D/2).$$
Hence for $x\in M$ and $r\ge D/2$, we get:
$$\int_{B(x,r)}q\le\int_M q\le \gamma\sup_{a\in A}\int_{B(x,D/2)}q.$$
And
$$\fint_{B(x,r)}q\le \gamma^2\sup_{a\in A}\fint_{B(x,D/2)}q.$$
Hence
\begin{equation*}
\int_{D/2}^\infty r e^{-\frac{r^2}{6R^2}}\left(\fint_{B(x,r)} q\right)dr\le \gamma^2 6R^2e^{-\frac{D^2}{24R^2}}\sup_{a\in A}\fint_{B(x,D/2)}q.
\end{equation*}
Similarly, 
\begin{equation*}\begin{split}
\int_{D/2}^D r e^{-\frac{r^2}{7R^2}}\left(\fint_{B(x,r)} q\right)dr&\ge \int_{D/2}^D r e^{-\frac{r^2}{7R^2}}dr\gamma^{-1}\left(\fint_{B(a,D/2)} q\right)\\
&\ge 7R^2 e^{-\frac{D^2}{28R^2}}\gamma^{-1}\left(\fint_{B(x,D/2)} q\right)
\end{split}\end{equation*}
If we define 
$$\Gamma:=\sup_x\int_0^D r e^{-\frac{r^2}{7R^2}}\left(\fint_{B(x,r)} q\right)dr$$
we easily get for any $x\in M$:
\begin{equation*}\begin{split}\int_0^\infty r e^{-\frac{r^2}{6R^2}}\left(\fint_{B(x,r)} q\right)dr&\le \Gamma+\gamma^26R^2e^{-\frac{D^2}{24R^2}}\sup_{a\in A}\fint_{B(a,D/2)}q\\
&\le \Gamma+\gamma^3e^{\frac{-D^2}{168R^2}}\Gamma.
\end{split}\end{equation*}
There is a constant $c$ (depending only on $C$) such that if $R\le D$ then
$$\gamma^3e^{\frac{-D^2}{168R^2}}\le  e^{\frac{-D^2}{168R^2}+3C\frac{D}{R}}\le c$$
Hence the result.
\endproof
\begin{thm}\label{katovolclosed} There is a constant $\upepsilon_n>0$ depending only on $n$ such that if  $(M^n,g)$ is a closed Riemannian manifold  such that some $R\le D$ and such that for all $x\in M$:
$$
\int_0^{D} r e^{\frac{-r^2}{7R^2}}\left(\fint_{B(x,r)} \ricm \right)dr\le \upepsilon_n ,$$
then $(M,g)$ satisfies 
 $\kappa_{R^2}(\ricm)\le \frac{1}{16n}$.
\end{thm}
\proof In \cite{Carron-16}, it has been shown that if 
 $\kappa_{R^2}(\ricm)\le \frac{1}{16n}$ then $(M^n,g)$ satisfies the condition (\ref{doubling}) and (\ref{heatkernel})  with constants that depend only on $n$: ($\nu=e^2n$).
 \begin{enumerate}[i)]
\item $\forall x\in M, 0\le r\le \rho\le R\colon \frac{\vol B(x, \rho)}{\vol B(x,r)}\le \uptheta_n\left(\frac{ \rho}{r}\right)^{\nu}$.
\item $\forall x,y\in M, t\in (0,R^2]\colon H(t,x,y)\le \frac{C_n}{\vol B(x,\sqrt{t})}e^{-\frac{d^2(x,y)}{5t}}.$
\end{enumerate}
Hence there is a constant $\uplambda_n$ depending only on $n$ such that 
$$\kappa_{R^2}(\ricm)\le \uplambda_n\sup_{x\in M}\int_0^{D} r e^{\frac{-r^2}{7R^2}}\left(\fint_{B(x,r)} \ricm \right)dr.$$
Hence if $R$ is the greatest real number such that 
$\kappa_{R^2}(\ricm)\le \frac{1}{16n}$ then we get $\kappa_{R^2}(\ricm)= \frac{1}{16n}$ and necessarily 
$$\frac{1}{16n}\le \uplambda_n\sup_{x\in M}\int_0^{D} r e^{\frac{-r^2}{7R^2}}\left(\fint_{B(x,r)} \ricm \right)dr.$$
By contraposition, we get that if 
$$\sup_{x\in M}\int_0^{D} r e^{\frac{-r^2}{7R^2}}\left(\fint_{B(x,r)} \ricm \right)dr\le \frac{1}{16 n \uplambda_n}$$ then
$$\kappa_{R^2}(\ricm)\le \frac{1}{16n}.$$
\endproof

From the above we can obtain many results from \cite{Carron-16} and \cite {Rose:2016aa}, we prefer to refrain from stating all of them and concentrate on few of them:
\begin{prop} There is a $\upeta_n$ such that if for all $x\in M$:
$$\int_0^{D} r \left(\fint_{B(x,r)} \ricm \right)dr\le \upeta_n$$
then the first Betti number of $M$ is smaller then $n$:
$$b_1(M)\le n.$$
\end{prop}

\begin{prop} Assume that  $(M^n,g)$ is a closed Riemannian manifold with diameter $D$ and introduce
$\mathbf{I}(M,g)$ defined by
$$\sup_{x\in M} \int_0^{D} r \left(\fint_{B(x,r)} \ricm^p \right)dr=\left(\mathbf{I}(M,g)\, D\right)^{2(p-1)}.$$ 
There is an explicit constant $\uptheta(n,p)$ such that for $R=\uptheta(n,p)D/\mathbf{I}(M,g)$:
$$\kappa_{R^2}(\ricm)\le \frac{1}{16n}$$ and with $\xi=\mathbf{I}(M,g)/\uptheta(n,p)$, one gets the isoperimetric inequality:
 for any open subset $\Omega\subset M$ with smooth boundary and $\vol \Omega\le \frac12 \vol M$,
we have
\begin{align*}
1\leq c_n^{1+\xi} D \vol(M) I_n(M).
\end{align*}

\end{prop} 
\proof With the H\"older inequality, one easily gets (with $q=p/(p-1)$):
\begin{align*}
\int_0^Dre^{-\frac{r^2}{7R^2}}\left(\fint_{B(x,r)}\ricm \right)dr&\le \int_0^Dre^{-\frac{r^2}{R^2}}\left(\fint_{B(x,r)}\ricm^p\right)^{\frac{ 1}{p}}dr\\
&\le \left(\int_0^Dre^{-q\frac{r^2}{7R^2}}dr\right)^{\frac{1}{ q}}\left(\mathbf{I}(M,g)\, D\right)^{\frac{2}{ q}}\\
&\le \left(\frac{7}{q}R^2\right)^{\frac 1 q}\left(\mathbf{I}(M,g)\, D\right)^{\frac{2}{ q}}.
\end{align*}
The first result follows now easily from the \tref{katovolclosed}. The second result follows from \cref{estikatoclosed} and \tref{isoperimetryprofile}. Indeed, 
we have 
$$\kappa_{R^2}(\ricm^p)\le \lambda_n \left(\mathbf{I}(M,g)\, D\right)^{2(p-1)}.$$
\endproof


\end{document}